\documentclass[final]{siamltex}
\usepackage{amsmath}
\usepackage{amssymb}
\usepackage[ruled,section]{algorithm}
\usepackage{algorithmic}

\usepackage[square,comma,numbers,sort&compress]{natbib} 

\usepackage{ifpdf}
\ifpdf 
\usepackage[pdftex]{graphicx}      
\DeclareGraphicsExtensions{.pdf}   
\usepackage[pdftex]{thumbpdf}      
\usepackage[pdftex,                
pdfstartview=FitH,%
bookmarks=true,
bookmarksnumbered=true,
bookmarksopen=true,%
hypertexnames=false,
breaklinks=true,
  colorlinks=true,%
  linkcolor=blue,anchorcolor=blue,%
  citecolor=blue,filecolor=blue,%
  menucolor=blue,pagecolor=blue]%
{hyperref} 
\hypersetup{
pdfauthor   = {Andrew V. Knyazev and Ilya Lashuk},
pdftitle    = {Steepest descent and conjugate gradient methods with variable preconditioning},
pdfsubject  = {Paper},
pdfkeywords = {Steepest descent, conjugate gradient, iterative method, inner-outer iterations, variable preconditioning, random preconditioning, preconditioner, condition number, linear systems, circular cone.},
}
\else 
\usepackage[dvips]{graphicx}       
\DeclareGraphicsExtensions{.eps}   
\usepackage[ps2pdf]{thumbpdf}      
\usepackage[ps2pdf,                
hypertexnames=false,
breaklinks=true,
  colorlinks=true,%
  linkcolor=blue,anchorcolor=blue,%
  citecolor=blue,filecolor=blue,%
  menucolor=blue,pagecolor=blue]%
{hyperref} 
\hypersetup{
pdfauthor   = {Andrew V. Knyazev and Ilya Lashuk},
pdftitle    = {Steepest descent and conjugate gradient methods with variable preconditioning},
pdfsubject  = {Paper},
pdfkeywords = {Steepest descent, conjugate gradient, iterative method, inner-outer iterations, variable preconditioning, random preconditioning, preconditioner, condition number, linear systems, circular cone.},
pdfcreator  = {LaTeX with hyperref package},
pdfproducer = {dvips + ps2pdf}
}
\fi 

\newcommand{\D}{\displaystyle}
\newcommand{\lb}{\left(}
\newcommand{\rb}{\right)}
\newcommand{\kmax}{\kappa_{\max}}
\newcommand{\vectornorm}[1]{\left\|#1\right\|}
\newcommand{\vnorm}{\vectornorm}

\title{Steepest descent and conjugate gradient methods with variable preconditioning
\thanks{%
Compiled \today. 
Received by the editors November 17, 2006;
accepted for publication (in revised form) ????????????;
published electronically ??????????.
Preliminary version of this paper is 
available as a technical report \cite[v1-v3]{kl06}. 
}}

\author{Andrew V. Knyazev\thanks{Department of
Mathematical Sciences University of Colorado at Denver
and Health Sciences Center,
P.O. Box 173364, Campus Box 170, Denver, CO 80217-3364
(Andrew.Knyazev[at]cudenver.edu, http://math.cudenver.edu/\~{}aknyazev/).
This material is based upon work supported by the National Science
Foundation awards DMS 0208773 and 0612751}
\and Ilya Lashuk\thanks{Department of
Mathematical Sciences University of Colorado at Denver
and Health Sciences Center,
P.O. Box 173364, Campus Box 170, Denver, CO 80217-3364 (Ilya.Lashuk[at]cudenver.edu).}}

\begin{document}
\setcounter{page}{1}
\maketitle
\begin{abstract}
We analyze the conjugate gradient (CG) method with variable preconditioning for solving a linear system with a real symmetric positive definite (SPD) matrix of coefficients $A$. We assume that the preconditioner is SPD on each step, and that the condition number of the preconditioned system matrix is bounded above by a constant 
independent of the step number. We show that the CG method with variable preconditioning under this assumption may not give improvement, compared to the steepest descent (SD) method. 
We describe the basic theory of CG methods with variable preconditioning with the emphasis on ``worst case'' scenarios, and provide complete proofs of all facts not available in the literature. We give a new elegant geometric proof of the SD convergence rate bound. 
Our numerical experiments, comparing the preconditioned SD and CG methods, not only support and illustrate our theoretical findings, but also reveal two surprising and potentially practically important effects. 
First, we analyze variable preconditioning in the form of inner-outer iterations. In previous such tests, the unpreconditioned CG inner iterations are applied to an artificial system with some fixed preconditioner as a matrix of coefficients. We test a different  scenario, where the unpreconditioned CG inner iterations solve linear systems with the original system matrix $A$. 
We demonstrate that the CG-SD inner-outer iterations perform as well as the CG-CG inner-outer iterations in these tests. 
Second, 
we compare the CG methods using a two-grid preconditioning with fixed and randomly chosen coarse grids, and observe that the fixed preconditioner method is twice as slow as the method with 
random preconditioning. 
\end{abstract}
\begin{AM}
65F10  
\end{AM}
\pagestyle{myheadings}
\thispagestyle{plain}
\markboth{ANDREW V. KNYAZEV AND ILYA LASHUK}
{SD AND CG METHODS WITH VARIABLE PRECONDITIONING}

\section{Introduction}

Preconditioning, a transformation,  usually implicit, of the original linear system aiming at accelerating the convergence of the approximations to the solution, is typically a necessary part of an efficient iterative technique. Modern preconditioning, e.g.,\ based on so-called algebraic multilevel and domain decomposition methods, attempts to become as close to a ``black box'' ideal  of direct solvers as possible. In this attempt, the mathematical structure of the preconditioner, which in the classical case is regarded as some linear transformation, may become very complex, in particular, the linearity can be easily lost, e.g., if the preconditioning itself involves ``inner'' iterative solvers. The fact that the preconditioner may be nonlinear, or variable, i.e.,\ changing from iteration to iteration, may drastically affect the known theory as well as the practical behavior of preconditioned iterative  methods and therefore needs special attention. Our main result is that the conjugate gradient (CG) method with variable preconditioning in certain situations may not give improvement, compared to the steepest descent (SD) method for solving a linear system with a real symmetric positive definite (SPD) matrix of coefficients. We assume that the preconditioner is SPD on each step, and that the condition number of the preconditioned system matrix is bounded above by a constant. 

Let us now introduce the notation, so that we can 
formulate the main result mathematically.  
Let $A$ be a real SPD matrix, 
$(x,y)$ be the standard inner product of real vectors $x$ and $y$, so that 
$(Ax,y)=(x,Ay)$, and let 
$\|x\|=\sqrt{(x,x)}$ be the corresponding vector norm.
We also use $\|\cdot\|$ to denote the operator norm. 
The $A$-inner product and the A-norm are denoted by 
$(x,y)_A=(x,Ay)$ and $\|x\|_A=\sqrt{(x,x)_A}.$   

We consider a family of iterative methods 
to obtain a sequence of approximate solutions $x_k$ of a linear system $Ax=b$ 
and use the $A$-norm to measure the error $e_k=x-x_k$. 
The SD and CG methods are well-known iterative procedures that fit into our framework. 
To accelerate the convergence of the error $e_k$ to zero we introduce preconditioning, i.e.,\ on every iteration $k$ an operator $B_k$, called the preconditioner, possibly different for each iteration $k$, is applied 
to the residual $r_k=b-Ax_k$. A general algorithm, 
which includes the preconditioned SD or CG (PSD or PCG respectively) methods 
as particular cases, can be presented as follows, 
e.g., \citet[p. 540]{axel_MR1276069} and \citet[Algorithm 5.3]{MR1269945}: 
given $A$, $b$, $\left\{B_k\right\}$, $\left\{m_k\right\}$, $x_0$, for $k=0,1,\ldots$:
$r_k=b-Ax_k, s_k=B_k^{-1} r_k,$ and 
\begin{equation}\label{m} 
\D p_k=s_k-\sum_{l=k-m_k}^{k-1}\frac{\left(As_k,p_{l}\right)}{\left(Ap_{l},p_{l}\right)}p_{l},\quad
\D x_{k+1}=x_k+\frac{\left(r_k,p_k\right)}{\left(Ap_k,p_k\right)}p_k,
\end{equation}
where 
\begin{equation}
\label{ass}
0\le m_k \le k \mbox{ and }  m_{k+1} \leq m_k+1. 
\end{equation}
The latter condition is highlighted in \citet[p. 1447, line 1]{notay_MR1797890} 
and ensures that the formula for $p_k$ in (\ref{m}) performs the standard 
Gram--Schmidt $A$-orthogonalizations to previous search directions, 
which are already pairwise  $A$-orthogonal. 
The full orthogonalization 
that performs explicit $A$-orthogonalizations to all previous search directions
corresponds to $m_k=k$. 
Choosing $m_k=\min\{k,1\}$ gives 
the PCG method, e.g.,\ described in \citet[IPCG Algorithm]{golub_ye_MR1740397}. 
The connection of this 
PCG method to the commonly used PCG algorithm 
is discussed in section \ref{practical_pcg_sect} 
following \citet[Remark 2.3]{golub_ye_MR1740397}. 
The shortest recursion $m_k=0$ leads to the standard PSD method.  

It is well-known, e.g.,\ \citet[p. 34]{MR1396083} and \citet[Sec. 11.1.2, p. 458]{axel_MR1276069}
that if the preconditioner is SPD and fixed, $B_k=B=B^{\ast}>0$, a preconditioned method, 
such as  (\ref{m}), 
using the preconditioner $B$ can be viewed as the corresponding unpreconditioned method applied to the preconditioned system $B^{-1}Ax=B^{-1}b$ in the $B$-based inner product ${( x,y)}_B={(x,By)}$. This implies that the theory obtained for unpreconditioned methods remains valid for preconditioned methods, in particular, 
the $A$-orthogonalization terms with $l<k-1$ in the sum in (\ref{m}) vanish in exact arithmetic, e.g.,\ 
\citet[Sec. 11.2.6, Th. 11.5]{axel_MR1276069}.
The situation changes dramatically, however, if different preconditioners $B_k$ are used in the PCG method. 

The present paper concerns the behavior of method (\ref{m}), where the preconditioner $B_k$ varies from step to step, but remains SPD on each step and the spectral condition number 
$
\D\kappa\lb B_k^{-1}A\rb={\lambda_{\max}\lb B_k^{-1}A\rb}/{\lambda_{\min}\lb B_k^{-1}A\rb}
$
is bounded above by some constant $\kmax$ independent of the step number $k$. We note that the matrix $B_k^{-1}A$ is SPD with respect to, e.g.,\ the $B_k$ inner product, so its eigenvalues are real positive.
Let us highlight that our assumption 
$\kappa\lb B_k^{-1}A\rb\leq\kmax$ can be equivalently written as  
$\|I-B_k^{-1}A\|_{B_k}\leq\gamma$ with $\kmax=(1+\gamma)/(1-\gamma),$ 
assuming without loss of generality that $B_k$ is scaled such that 
${\lambda_{\max}\lb B_k^{-1}A\rb}+{\lambda_{\min}\lb B_k^{-1}A\rb}=2.$
Here, we only deal with methods that are 
invariant with respect to scaling of $B_k$. 

The main result of this paper is that the preconditioned method (\ref{m}) with (\ref{ass}) turns into the PSD method with the worst possible convergence rate on every iteration, if the preconditioners $B_k$ satisfying our assumption $\kappa\lb B_k^{-1}A\rb\leq\kmax$ are chosen in a special way. 
We explicitly construct a variable preconditioner that slows down the CG method to the point that the worst linear convergence rate of the SD method is recovered. 
Thus one can only guarantee that the convergence rate for the method (\ref{m}) with (\ref{ass}) 
is just the same as for the PSD method, $m_k=0$, obtained 
in \citet{MR0023126} and reproduced, e.g.,\ in 
\citet[Ch. XV]{Kantorovich-1964}: 
\begin{equation}
\label{classic_psd_bound}
\D\frac{\vnorm{e_{k+1}}_A}{\vectornorm{e_k}_A}\le\frac{\kappa_{\max}-1}{\kappa_{\max}+1}. 
\end{equation}
Our proof is geometric and is based on the simple fact, proved in  section \ref{suppl_sect}, 
that a nonzero vector multiplied by all SPD matrices with a condition number bounded by a constant generates a pointed circular cone. 
We apply this fact on every iteration to the current residual vector, 
which becomes the center of the cone, so all points in the cone 
correspond to all possible preconditioned residuals. 
In a somewhat similar way, \citet{golub_ye_MR1740397}
use the angle between the exact and the 
perturbed preconditioned residuals. 
In the CG method context, this cone has a nontrivial intersection with the subspace $A$-orthogonal to all previous search directions. So on each iteration we can choose a preconditioner with the {a priori} chosen quality, determined by $\kappa_{\max},$ that makes enforcing $A$-orthogonality with respect to all previous search directions useless. 
 
Basic properties of method (\ref{m}), most importantly the local optimality, are derived in section \ref{prop_sect}.   
In section \ref{alter_sd_sect} we apply our results from section \ref{suppl_sect} about the cone to obtain a new proof of estimate (\ref{classic_psd_bound}). 
In section \ref{main_sect} we analyze the convergence of the PCG method with variable preconditioning and prove our main result. 
We assume real arithmetic everywhere in the paper, except for 
section \ref{complex_case_sect}, where we show that our main results also hold for complex Hermitian positive definite matrices. 
In section \ref{practical_pcg_sect} we consider two particular PCG algorithms that are often used in practice and describe their behavior with variable preconditioning. 

Our numerical experiments in section \ref{num_exp_sect} comparing the preconditioned SD and CG methods support and illustrate our theoretical findings, and also reveal some potentially practically important effects. In subsection \ref{num_exp_ssect1}, we test the widely used modification of the CG method with a simplified formula for the scalar $\beta_k$ from section \ref{practical_pcg_sect} and demonstrate that variable preconditioning can make this modification much slower than even the SD method. 
In subsection \ref{num_exp_ssect2}, we analyze inner-outer iterations as variable preconditioning. 
Finally, in subsection \ref{num_exp_ssect3}, we demonstrate that variable preconditioning may surprisingly accelerate the SD and the CG compared to the use of fixed preconditioning in the same methods. 
 

Different aspects of variable preconditioning are considered, e.g.,\ 
in \citet{axel_vass_MR1121697,MR1269945,axel_MR1276069}, where rather general nonlinear preconditioning is introduced, and 
in \citet{notay_MR1797890,golub_ye_MR1740397} that mainly deal with the case when the preconditioner on each iteration approximates a fixed operator. In \citet{axel_vass_MR1121697,axel_MR1276069,notay_MR1797890,golub_ye_MR1740397}, convergence estimates for some iterative methods with variable preconditioning are proved. For recent results and other aspects of variable preconditioning see \citet{MR1974182,MR2058070,MR2179897} and references there. No attempts are apparently made in the literature to obtain a result similar to ours, even though it should appear quite natural and somewhat expected to experts in the area, after reading this paper. 

\section{Pointed circular cones represent sets of SPD matrices with varying condition numbers}
\label{suppl_sect}
For a pair of real non-zero 
vectors $x$ and $y$ we define the angle between $x$ and $y$ in the usual way as 
$$\D\angle(x,y)=\arccos\left(\frac
{(x,y)} {\vectornorm{x}\vectornorm{y}}\right)
\in[0,\pi].$$ 
The following theorem is inspired by \citet[Lemma 2.3]{neymeyr_MR1804114}.
\begin{theorem} \label{cone_th}
The set $\{Cx\}$, where $x$ is a fixed nonzero real vector and $C$ runs through all SPD matrices with 
condition number $\kappa(C)$ bounded above by some $\kappa_{\max}$, 
is a pointed circular cone, specifically, 
$$
\{Cx: \quad C=C^*>0,\, \kappa(C)\leq  \kappa_{\max} \}= \left\{ 
y: \quad \sin \angle\lb x,y\rb\leq\frac{\kappa_{\max}-1}{\kappa_{\max}+1}
\right\}.
$$
\end{theorem}

Theorem \ref{cone_th} can be proved by constructing our cone as the smallest 
pointed cone that includes the ball considered in \citet[Lemma 2.3]{neymeyr_MR1804114}.   
Preparing for section \ref{complex_case_sect}
that deals with the complex case, not covered in \citet{neymeyr_MR1804114}, 
we provide a direct proof here based on  the following two lemmas. 
The first lemma is simple and states that the set in question cannot be larger than the cone: 
\begin{lemma}
\label{sd-like_lemma}
Let $x$ be a non-zero real vector, let $C$ be an SPD matrix with spectral condition number $\kappa\lb C\rb$. Then 
$
\sin \angle\lb x,Cx\rb\leq
{(\kappa\lb C\rb-1)}/{(\kappa\lb C\rb+1)}.
$
\end{lemma}
\begin{proof}
Denote $y=Cx$. 
We have $\lb x,Cx\rb = \lb y,C^{-1}y\rb>0$ since $C$ is SPD, so $y\neq 0$ and $\angle\lb x,y\rb<{\pi}/{2}$.
A positive scaling of $C$ and thus of $y$ is obviously irrelevant, so let us 
choose $y$ to be the orthogonal projection of $x$ onto the 1-dimensional subspace spanned by the original $y$. 
Then from elementary 2D geometry it follows that $\vnorm{y-x}=\vnorm{x}\sin\angle\lb x,y\rb$. 
The orthogonal projection of a vector onto a subspace is the best approximation to the vector from the subspace, 
thus  
$$\vnorm{x}\sin\angle\lb x,y\rb=\vnorm{y-x}\le\vnorm{sy-x}=\vnorm{sCx-x}\le\vnorm{sC-I}\vnorm{x}$$ 
for any scalar $s$, where 
$I$ is the identity. Taking
$s=2/\lb{\lambda_{\max}\lb C\rb+\lambda_{\min}\lb C\rb}\rb,$
where $\lambda_{\min}\lb C\rb$ and $\lambda_{\max}\lb C\rb$ are the minimal and maximal eigenvalues of $C$, respectively, 
we get 
$\vnorm{sC-I} = (\kappa\lb C\rb-1)/(\kappa\lb C\rb+1)$. 
\end{proof}
 
The second lemma implies that every point in the cone can be represented as $Cx$ for some SPD matrix $C$ with $\kappa(C)$ determined by the opening angle of the cone. 
\begin{lemma}
\label{neymeyr-like_lemma}
Let $x$ and $y$ be non-zero real vectors, such that $\angle\lb x,y\rb\in\left[0,\frac{\pi}{2}\right)$. Then there exists an SPD matrix $C$, such that $Cx=y$ and 
$$
\frac{\kappa\lb C\rb-1}{\kappa\lb C\rb+1}=\sin\angle\lb x,y\rb. 
$$
\end{lemma}
\begin{proof}
Denote $\alpha=\angle\lb x,y\rb$. 
A positive scaling of vector $y$ is irrelevant, so as in the previous proof we  
choose $y$ to be the orthogonal projection of $x$ onto the 1-dimensional subspace spanned by the original $y$, then  
$\vnorm{y-x}=(\sin\alpha) \vnorm{x}$, so the vectors $y-x$ and $(\sin\alpha) x$ are of the same length. 
This implies that there exists a Householder reflection $H$ such that $H\lb\lb\sin\alpha\rb x\rb= y-x$, cf. \citet[Lemma 2.3]{neymeyr_MR1804114}, 
so $\lb I+\lb\sin\alpha\rb H\rb x= y$. We define $C=I+\lb\sin\alpha\rb H$ to get $Cx=y.$ 
Any Householder reflection is symmetric and has only two distinct eigenvalues $\pm 1$, so 
$C$ is also symmetric and has only two distinct positive eigenvalues $1\pm\sin\alpha$, 
as $\alpha\in\left[0,{\pi}/{2}\right)$, and we conclude that  $C>0$ and 
$\kappa\lb C\rb=(1+\sin\alpha)/(1-\sin\alpha)$.
\end{proof}

\section{Local optimality of the method with  variable preconditioning}
\label{prop_sect}
Here we discuss some basic properties of method (\ref{m}) with (\ref{ass}). 
We derive a simple, but very useful, error propagation identity in Lemma \ref{error_transition_lemma}. 
We prove in Lemma \ref{local_ortho} that the method is well-defined and has a certain local $A$-orthogonality property, formulated without a proof in  
\citet[formulas (2.1)-(2.2)]{notay_MR1797890} and in the important 
particular case  $m_k=\min\{k,1\}$ proved in 
\citet[Lemma 2.1]{golub_ye_MR1740397}. 
Using the local $A$-orthogonality property of Lemma \ref{local_ortho}, 
we prove the local A-optimality property in Lemma \ref{lo_l} 
by generalizing the result of \citet[Proposition 2.2]{golub_ye_MR1740397}.
Finally, we derive a trivial Corollary \ref{lo_csd} from Lemma \ref{lo_l}, 
which uses the idea from \citet[p. 1309]{golub_ye_MR1740397} of comparison with the 
PSD method, $m_k=0$. 

The material of this section is inspired by \citet{golub_ye_MR1740397}
and may be known to experts in the field, e.g.,\ some even more general facts can be found in 
\citet[Sec. 12.3.2, Lemma 12.22]{axel_MR1276069}. 
We provide straightforward and complete proofs here suitable for a general audience.
\begin{lemma}
\label{error_transition_lemma}
Let $A$ and $\{B_k\}$ be SPD matrices. Suppose $p_k$ in method (\ref{m}) is well-defined and nonzero. Then 
\begin{equation}
\label{error_transition}
\D e_{k+1}=e_k-\frac{\left(Ae_k,p_k\right)}{\left(Ap_k,p_k\right)}p_k.
\end{equation}
\end{lemma}
\begin{proof}
Recall that $e_k=A^{-1}b-x_k$ and thus $r_k=Ae_k$. Then (\ref{error_transition}) 
follows immediately from the last formula in (\ref{m}).
\end{proof}	
\begin{lemma}
\label{local_ortho}
Let $A$ and $\{B_k\}$ be SPD matrices and $\{m_k\}$ satisfies (\ref{ass}). Then the error, the preconditioned residual, and the direction vectors generated by method (\ref{m}) before the exact solution is obtained are well-defined and satisfy
\begin{equation}
\label{dirs_Aorth}
\lb p_i,p_j\rb_A=0,\;k-m_k\le i<j\le k,
\end{equation}
\begin{equation}
\label{A_ortho_prop}
\lb e_{k+1},s_k\rb_A=\lb e_{k+1},p_i\rb_A=0,\; k-m_k\le i\le k. 
\end{equation}
\end{lemma}
\begin{proof}
We first notice that (\ref{error_transition}) for any $k$ obviously implies
\begin{equation}
\label{error_dir_orth}
\lb e_{k+1},p_k\rb_A=0.
\end{equation}
For the rest of the proof we use an induction in $k$. 
Let us take $k=0$ and suppose $x_0\ne x$, then $r_0\ne 0$ and $s_0\ne 0$
since $B_0$ is SPD. By (\ref{ass}), $m_0=0$ and thus $p_0=s_0 \ne 0$,  
so in the formula for $x_{k+1}$ we do not divide by zero, i.e.,\ 
$x_{k+1}$ is well defined. 
There is nothing to prove in (\ref{dirs_Aorth}) for $k=0$ since $m_0=0$.
Formula (\ref{error_dir_orth}) implies $\lb e_{1},p_0\rb_A=\lb e_{1},s_0\rb_A=0$, 
 i.e.,\ (\ref{A_ortho_prop}) holds for $k=0$. 
This provides the basis for the induction.

Suppose the statement of the lemma holds for $k-1$, which is the induction hypothesis,
i.e.,\ up to the index $k-1$ all quantities are well defined and 
\begin{equation}
\label{prev_dirs}
\lb p_i,p_j\rb_A=0,\;k-1-m_{k-1}\le i<j\le k-1,
\end{equation}
\begin{equation}
\label{induct_ortho}
\lb e_k,s_{k-1}\rb_A=\lb e_k,p_{i}\rb_A=0,\; k-1-m_{k-1}\le i\le k-1. 
\end{equation}

We now show by contradiction that $x_k\ne x$ implies $p_k\ne 0$. Indeed. if
$p_k=0$ then $s_k$ is a linear combination of $\D p_{k-m_k},\ldots,p_{k-1}$. However, since $m_k\le m_{k-1}+1$,  it follows from (\ref{induct_ortho}) that
\begin{equation}
\label{almost_ortho}
\lb e_k,p_{i}\rb_A=0,\; k-m_{k}\le i\le k-1. 
\end{equation}
Then we have $\lb s_k,e_k\rb_A=0$. At the same time, since the matrix $B_k^{-1}A$ is $A$-SPD, $s_k=B_k^{-1}Ae_k$ cannot be $A$-orthogonal to $e_k$ unless $s_k=e_k=0$, i.e.,\ $x_k=x$. 

Next, we prove (\ref{dirs_Aorth}) by showing that 
the formula for $p_k$ in (\ref{m}) is a valid step of 
the Gram-Schmidt orthogonalization process with respect to the $A$-based inner product. 
If $m_k=0,$ there is nothing to prove. If $m_k=1$ then (\ref{dirs_Aorth}) gets reduced to $\lb p_k,p_{k-1}\rb_A=0,$ which follows from the formula for $p_k$ in (\ref{m}). If $m_k\ge 2$ then condition (\ref{ass}) implies that vectors $\D p_{k-m_k},\ldots,p_{k-1}$
are among the vectors
$\D p_{k-1-m_{k-1}},\ldots,p_{k-1}$ and therefore are already $A$-orthogonal by 
the induction assumption (\ref{prev_dirs}). Then the formula for $p_k$ in (\ref{m}) is indeed a valid step of the Gram-Schmidt orthogonalization process with respect to the $A$-based inner product, so (\ref{dirs_Aorth}) holds.

It remains to prove (\ref{A_ortho_prop}). We have already established 
(\ref{dirs_Aorth}), and (\ref{error_dir_orth})--(\ref{almost_ortho}).
Equalities (\ref{dirs_Aorth}) and (\ref{almost_ortho}) imply that $p_k$ and $e_k$ are $A$-orthogonal to $\D p_{k-m_k},\ldots,p_{k-1}$. Equality (\ref{error_transition}) implies that $e_{k+1}$ is a linear combination of $e_k$ and $p_k$. Thus, we have
$
\lb e_{k+1},p_{i}\rb_A=0,\; k-m_{k}\le i\le k-1.
$
Finally, it is enough to notice that $s_k$ is a linear combination of $p_k,p_{k-1},\ldots,p_{k-m_k}$, so $(e_{k+1},s_k)_A=0$.
\end{proof}

We now use Lemma \ref{local_ortho} to prove the local optimality of method 
(\ref{m}) with (\ref{ass}), which generalizes 
the statement of \citet[Proposition 2.2]{golub_ye_MR1740397}.
\begin{lemma}\label{lo_l}
Under the assumptions of Lemma \ref{local_ortho}, 
$$
\vnorm{e_{k+1}}_A=\min_{p\in\text{span}\left\{s_k,p_{k-m_k},\ldots,p_{k-1}\right\}}{\vnorm{e_k-p}_A}.
$$
\end{lemma}
\begin{proof}
We get $e_{k+1}\in e_k+\text{span}\left\{s_k,p_{k-m_k},\ldots,p_{k-1}\right\}$
from the formula for $p_k$ in (\ref{m}) and (\ref{error_transition}). 
Putting this together with $A$-orthogonality relations (\ref{A_ortho_prop}) 
of the vector $e_{k+1}$ with all vectors that span the subspace finishes the proof.   
\end{proof}

Two important corollaries follow immediately from Lemma \ref{lo_l} by analogy with 
\citet[Proposition 2.2]{golub_ye_MR1740397}.
\begin{corollary}\label{lo_csd}
The $A$-norm of the error $\vectornorm{e_{k+1}}_A$ in method (\ref{m}) with (\ref{ass}) is bounded 
above by the $A$-norm of the error of one step of the PSD method, $m_k=0$, using the same $x_k$ 
as the initial guess 
and $B_k$ as the preconditioner, i.e., specifically,  
$\vnorm{e_{k+1}}_A\leq\min_{\alpha}{\vnorm{e_k-\alpha s_k}_A}.$
\end{corollary}
\begin{corollary}\label{lo_chb}
Let $m_k>0$, then 
the $A$-norm of the error $\vectornorm{e_{k+1}}_A$ in method (\ref{m}) with (\ref{ass}) 
for $k>0$ satisfies
$\vnorm{e_{k+1}}_A\leq\min_{\alpha, \beta}{\vnorm{e_k-\alpha s_k -\beta (e_k-e_{k-1})}_A}.$
\end{corollary}
\begin{proof}
Under the lemma assumptions, the formula for $p_k$ in (\ref{m}) and (\ref{error_transition}) imply that 
$
e_{k+1}\in e_k+\text{span}\left\{s_k,p_{k-1}\right\} = 
 e_k+\text{span}\left\{s_k,e_k-e_{k-1}\right\}, 
$
and the $A$-orthogonality relations (\ref{A_ortho_prop}) 
turn into $(e_{k+1},s_k)_A=0$ and  
$(e_{k+1},p_{k-1})_A = (e_{k+1},e_k-e_{k-1})_A=0$, 
so the vector $e_{k+1}$ is $A$-orthogonal to both vectors that span the subspace. 
As in the proof of Lemma \ref{lo_l}, 
the local  $A$-orthogonality implies the local $A$-optimality.
\end{proof}

Corollary \ref{lo_csd} allows us in section \ref{alter_sd_sect} to estimate the convergence rate 
of method (\ref{m}) with (\ref{ass}) by comparison with the 
PSD method, $m_k=0$,---this idea is borrowed from 
\citet[p. 1309]{golub_ye_MR1740397}. 
The results of Lemma \ref{lo_l} and Corollary \ref{lo_chb} seem to indicate that an improved  
convergence rate bound of method (\ref{m}) with (\ref{ass}) can be obtained, compared to the PSD method 
convergence rate bound that follows from Corollary \ref{lo_csd}. 
Our original intent has been to combine Corollary \ref{lo_chb} with
convergence rate bounds of the heavy ball method, in order to attempt
to prove such an improved convergence rate bound. 
However, our results of section \ref{main_sect} demonstrate that this improvement is impossible
under our only assumption $\kappa\lb B_k^{-1}A\rb\leq\kmax$, since one can construct 
such preconditioners $B_k$ that, e.g.,\ the minimizing value of $\beta$ in Corollary \ref{lo_chb} is zero, so 
Corollary \ref{lo_chb} gives no improvement compared to Corollary \ref{lo_csd}. 

\section{Convergence rate bounds for variable preconditioning}
\label{alter_sd_sect}
The classical \citet[Ch. XV]{Kantorovich-1964} convergence rate bound (\ref{classic_psd_bound}) for the PSD method
is ``local'' in the sense that it relates the $A$-norm of the error 
on two subsequent iterations and does not depend on previous iterations. 
Thus, it remains valid when the preconditioner $B_k$ changes from iteration to iteration, 
while the condition number $\kappa\lb B_k^{-1}A\rb$ is bounded above by some constant $\kmax$ independent of $k$. 
The goal of this section is to give an apparently new simple proof of the estimate (\ref{classic_psd_bound}) for the PSD method, 
based on our cone Theorem \ref{cone_th}, and to extend this statement to cover the general method (\ref{m}) with (\ref{ass}), using Corollary \ref{lo_csd}.  

We denote the angle between two real nonzero vectors with respect to the $A$-based inner product by 
\[
\D\angle_A\lb x,y\rb=\arccos\left(\frac{(x,y)_A}{\vectornorm{x}_A\vectornorm{y}_A}\right) \in [0,\pi]
\]
and express the error reduction ratio for the PSD method in terms of the angle with respect to the $A$-based inner product:  
\begin{lemma}
\label{sd_reduct_lemma}
On every step of the PSD algorithm, (\ref{m}) with $m_k=0$, the error reduction factor takes the form  
$
{\vectornorm{e_{k+1}}_A}/{\vectornorm{e_k}_A} = 
\sin(\angle_A(e_k,B_k^{-1}Ae_k)).$ 
\end{lemma}
\begin{proof}
By (\ref{A_ortho_prop}), we have $(e_{k+1},p_k)_A=0$ . Now, for  $m_k=0$, in addition, 
$p_k=s_k$, so  $0=(e_{k+1},p_k)_A=(e_{k+1},s_k)_A=(e_{k+1},x_{k+1}-x_k)_A$, i.e.,\ 
the triangle with vertices $x$, $x_k$, $x_{k+1}$ is right-angled in the $A$-inner product, 
where the hypotenuse is $e_k=x-x_k$. Therefore,  
$\vectornorm{e_{k+1}}_A / \vectornorm{e_{k}}_A = \sin(\angle_A(e_k,x_{k+1}-x_k)) = \sin(\angle_A(e_k,s_k))$, 
where $s_k=B_k^{-1}\lb b-Ax_k\rb=B_k^{-1}Ae_k$ by (\ref{m}). 
\end{proof}

Let us highlight that Lemma \ref{sd_reduct_lemma} provides an exact expression for the error reduction factor, not just a bound---we need this in the proof of Theorem \ref{main_theorem}
in the next section.
Combining the results of Lemmas \ref{sd-like_lemma} and \ref{sd_reduct_lemma} together immediately leads to 
(\ref{classic_psd_bound}) for the PSD method, where $m_k=0$. Finally, taking into account 
Corollary \ref{lo_csd}, by analogy with the arguments of \citet[p. 1309]{golub_ye_MR1740397} 
and decrypting a hidden statement in \citet[Lemma 3.5]{golub_ye_MR1740397}, 
we get 
\begin{theorem}\label{cr_th}
Convergence rate bound (\ref{classic_psd_bound}) 
holds for method (\ref{m}) with (\ref{ass}).
\end{theorem}

\section{The convergence rate bound is sharp}
\label{main_sect}
Here we formulate and prove the main result of the paper that 
one can only guarantee the convergence rate described by (\ref{classic_psd_bound})
for method (\ref{m}) with (\ref{ass}) with variable preconditioning 
if one only assumes $\kappa\lb B_k^{-1}A\rb\leq\kmax$.  
Let us remind the reader that (\ref{classic_psd_bound}) 
also describes the convergence rate for the PSD method, (\ref{m}) with $m_k=0$. 
We now show that adding more vectors 
to the PSD iterative recurrence results in no improvement in convergence, 
if a specially constructed set of variable preconditioners is used.  
\begin{theorem}
\label{main_theorem}
Let an SPD matrix $A$, vectors $b$ and $x_0$, and $\kappa_{\max}>1$ be given.
Assuming that the matrix size is larger than the number of iterations, one can choose a sequence of SPD preconditioners $B_k$, satisfying $\kappa(B_k^{-1}A)\le \kappa_{\max}$, such that 
method (\ref{m}) with (\ref{ass})
turns into the PSD method, (\ref{m}) with $m_k=0$, and on every iteration
\begin{equation}\label{e_eq}
\frac{\vectornorm{e_{k+1}}_A}{\vectornorm{e_k}_A}=\frac{\kappa_{\max}-1}{\kappa_{\max}+1}.
\end{equation}
\end{theorem}
\begin{proof}
We construct the sequence $B_k$ by induction. First, we choose any vector $q_0$, such that 
$\sin \angle_A(q_0,e_0)=(\kappa_{\max}-1)/(\kappa_{\max}+1).$
According to Lemma \ref{neymeyr-like_lemma} applied in the $A$-inner product, 
there exists an $A$-SPD matrix $C_0$ with condition number $\kappa(C_0)=\kappa_{\max}$, such that $C_0e_0=q_0$. 
We define the SPD $B_0=AC_0^{-1}$, then $\kappa(B_0^{-1}A)= \kappa(C_0)=\kappa_{\max}$. 
We have $s_k=B_k^{-1}Ae_k$, so such a choice of $B_0$ implies $s_0=q_0$. 
Also, we have $p_0=s_0$, i.e.,\ the first step is always a PSD step, thus, by Lemma \ref{sd_reduct_lemma} we have 
proved (\ref{e_eq}) for $k=0$. 
Note that $(e_{1},p_0)_A=0$ by (\ref{A_ortho_prop}).

Second, we make the induction assumption: let preconditioners $B_l$ for $l\le k-1$ be constructed, such that 
$
{\vectornorm{e_{l+1}}_A}/{\vectornorm{e_l}_A}=
{(\kappa_{\max}-1)}/{(\kappa_{\max}+1)}
$
and $(e_k, p_{l})_A=0$ hold for all $l\le k-1$.
The dimension of the space is greater than the total number of iterations by our assumption, so there exists a vector $u_k$, such that 
$(u_k,p_{l})_A=0$ for $l\le k-1$ and $u_k$ and $e_k$ are linearly independent. 
Then the 2D subspace spanned by $u_k$ and $e_k$ is $A$-orthogonal to $p_{l}$ for $l\le k-1$. 

Let us consider the boundary of the pointed circular cone 
made of vectors $q_k$ satisfying the condition 
$\sin \angle_A(q_k,e_k)=(\kappa_{\max}-1)/(\kappa_{\max}+1).$
This conical surface has a nontrivial intersection with the 
2D subspace spanned by $u_k$ and $e_k$, since $e_k$ is the cone axis.  
Let us choose vector $q_k$ in the intersection, 
This vector will be obviously $A$-orthogonal to $\D p_{l}$, $l\le k-1$. 

Applying the same reasoning as for constructing $B_0$, 
we deduce that there exists an SPD $B_k$ such that $\kappa(B_k^{-1}A)\le\kappa_{\max}$ and $B_k^{-1}Ae_k=q_k$. 
With such a choice of $B_k$ we have $s_k=q_k$. Since $q_k=s_k$ is $A$-orthogonal to $p_{l}$ for all $l\le k-1$, it turns out that $p_k=s_k$, 
no matter how $\{m_k\}$ are chosen. This means that $x_{k+1}$ is obtained from $x_k$ by a steepest descent step. 
Then we apply Lemma \ref{sd_reduct_lemma} and conclude that (\ref{e_eq}) holds.
We note, that $(e_{k+1},p_l)_A=0$ for all $l\le k$. Indeed, $(e_{k+1},p_l)_A=0$ for all $l\le k-1$ 
since $e_{k+1}$ is a linear combination of $e_k$ and $p_k=s_k=q_k$, both $A$-orthogonal to $p_l$ for $l\le k-1$. 
Finally, $(e_{k+1},p_k)_A=0$ by (\ref{A_ortho_prop}).
This completes the construction of $\{B_k\}$ by induction and thus the proof.
\end{proof}

Let us highlight that 
the statement of Theorem \ref{main_theorem} consists of two parts: 
first, it is possible to have the PCG method with variable preconditioning 
that converges not any faster than the PSD method with the same preconditioning; 
and second, moreover, it is possible that the PCG method with variable preconditioning
converges not any faster than the worst possible theoretical 
convergence rate for the PSD method described by (\ref{classic_psd_bound}).  
Numerical tests in section \ref{num_exp_sect} show that the former 
possibility is more likely than the latter. Specifically, 
we demonstrate numerically in subsection \ref{num_exp_ssect3} 
that the PCG and PSD methods with random 
preconditioning converge with the same speed, but both are much faster than 
what bound (\ref{classic_psd_bound}) predicts.  

\section{Complex Hermitian case}
\label{complex_case_sect}
In all other sections of this paper we assume for simplicity that matrices and vectors are real. However, our main results also hold when matrices $A$ and $\{B_k\}$ are complex Hermitian positive definite. In this section we discuss necessary modifications to statements and proofs in sections \ref{suppl_sect}, \ref{alter_sd_sect} and \ref{main_sect} in order to cover the complex Hermitian case.

In section \ref{suppl_sect}, the first thing to be changed is the definition of the angle between two non-zero vectors $x,y\in{\mathbb{C}}^n$, where  an absolute value is now taken,
\[
\angle(x,y)=\arccos\left|\frac{(x,y)}{\vectornorm{x}\vectornorm{y}}\right| \in [0,\frac{\pi}{2}],
\]
that makes the angle acute and invariant with respect to complex nonzero scaling of the vectors. 
Lemma \ref{sd-like_lemma} remains valid in the complex case:
\begin{lemma}
\label{sd-like_lemma_complex}
Let $x$ be a non-zero complex vector, and $C$ be a complex Hermitian positive definite matrix with the spectral condition number $\kappa\lb C\rb$, then 
$
\sin \angle\lb x,Cx\rb\leq
{(\kappa\lb C\rb-1)}/{(\kappa\lb C\rb+1)}.
$
\end{lemma}
\begin{proof}
Denote $y=Cx$ and 
let $\gamma=\lb y,x\rb/\vnorm{y}^2$ then $\gamma y$ is the projection of $x$ onto  
${\rm span}\{y\},$ and $\angle\lb x,y\rb= \angle\lb x,\gamma y\rb$. 
Moreover, $(x, \gamma y) = (x, y) \gamma = (x, y)(y,x)/\vnorm{y}^2$ is 
real---we need this fact later in the proof of Lemma \ref{neymeyr-like_lemma_complex}. 
We redefine $y$ to $\gamma y$. 
The rest of proof is exactly the same as that of Lemma \ref{sd-like_lemma}, since the identity $\vnorm{y-x}=\vnorm{x}\sin\angle\lb x,y\rb,$ where $y$ is scaled 
by a complex scalar to be the orthogonal  projection of $x$ onto ${\rm span}\{y\},$ holds in the complex case with the new definition of the angle. 
\end{proof}

Lemma \ref{neymeyr-like_lemma} and, thus, Theorem \ref{cone_th} 
do not hold in the complex case after the straightforward reformulation. A trivial counterexample is a pair of vectors $x\neq 0$ and $y=ix$---the angle between $x$ and $y$ is obviously zero, yet it is impossible that $y=Cx$ for any complex Hermitian matrix $C,$ since the inner product $\lb x,y\rb=-i{\vnorm{x}}^{2}$ is not a real number.
This counterexample also gives an idea for a simple fix: 
\begin{lemma}
\label{neymeyr-like_lemma_complex}
Let $x$ and $y$ be non-zero complex vectors, such that $\angle\lb x,y\rb\neq{\pi}/{2}$. Then there exists a complex Hermitian positive definite matrix $C$ and a complex scalar $\gamma$, such that $Cx=\gamma y$ and 
$
{(\kappa\lb C\rb-1)}/{(\kappa\lb C\rb+1)}=\sin\angle\lb x,y\rb. 
$
\end{lemma}
\begin{proof}
We first scale the complex vector $y$ as in the proof 
of Lemma \ref{sd-like_lemma_complex} to make 
$y$ to be the projection of $x$ onto ${\rm span}\{y\}.$
The rest of the proof is 
similar to that of Lemma \ref{neymeyr-like_lemma}, but we have to  
be careful working with the Householder reflection in the complex case, so we provide the complete proof. 

The redefined $y$ is the projection of $x$ onto ${\rm span}\{y\},$ thus,  
$\vnorm{y-x}=(\sin\alpha) \vnorm{x}$, so the vectors $u=y-x$ and $v=(\sin\alpha) x$ are of the same length. Moreover, their inner product $(u,v)$ 
is real, since $(x, y)$ is real, see the proof of Lemma \ref{sd-like_lemma_complex}.  
This implies that the Householder reflection 
$Hz = z - 2 (w,z) w$, where $w=(u-v)/\|u-v\|$,  
acts on $z=u$ such that $Hu=v$, i.e.,\ $H\lb\lb\sin\alpha\rb x\rb= y-x$, 
so $\lb I+\lb\sin\alpha\rb H\rb x= y$. We define $C=I+\lb\sin\alpha\rb H$ to get $Cx=y.$ 

The Householder reflection $H$ is Hermitian and has only two distinct eigenvalues $\pm 1$, so 
$C$ is also Hermitian and has only two distinct positive eigenvalues $1\pm\sin\alpha$, 
as $\alpha\in\left[0,{\pi}/{2}\right)$, and we conclude that  $C>0$ and 
$\kappa\lb C\rb=(1+\sin\alpha)/(1-\sin\alpha)$.
\end{proof}

The same change then makes Theorem \ref{cone_th} work in the complex case: 
\begin{theorem} \label{cone_th_c}
The set $\{\gamma Cx\}$, 
where $x$ is a fixed nonzero complex vector, 
$\gamma$ runs through all nonzero complex scalars, and $C$ runs through all 
complex Hermitian positive definite matrices with 
condition number $\kappa(C)$ bounded above by some $\kappa_{\max}$, 
is a pointed circular cone, specifically, 
$$
\{\gamma Cx: \quad \gamma\neq 0, C=C^*>0,\, \kappa(C)\leq  \kappa_{\max} \}= \left\{ 
y: \quad \sin \angle\lb x,y\rb\leq\frac{\kappa_{\max}-1}{\kappa_{\max}+1}
\right\}.
$$
\end{theorem}

Section \ref{prop_sect} requires no changes other then replacing ``SPD'' with 
``Hermitian positive definite.'' 
In section \ref{alter_sd_sect} we just   
change the definition of the $A$-angle to
\[
\angle_A\lb x,y\rb=\arccos\left|\frac{(x,y)_A}{\vectornorm{x}_A\vectornorm{y}_A}\right| \in [0,\frac{\pi}{2}], 
\]
and then Lemma \ref{sd_reduct_lemma} holds without any further changes.

Finally, the statement of Theorem \ref{main_theorem} from section \ref{main_sect} 
allows for a straightforward generalization:
\begin{theorem}
\label{main_theorem_complex}
Let a Hermitian positive definite matrix $A$, complex vectors $b$ and $x_0$, and $\kappa_{\max}>1$
be given. Assuming that the matrix size is larger than the number of iterations, one can choose a sequence of Hermitian positive definite preconditioners $B_k$, satisfying $\kappa(B_k^{-1}A)\le \kappa_{\max}$, 
such that method (\ref{m}) with (\ref{ass})
turns into the PSD method, (\ref{m}) with $m_k=0$, and on every iteration
\begin{equation}\label{e_eq_complex}
\frac{\vectornorm{e_{k+1}}_A}{\vectornorm{e_k}_A}=\frac{\kappa_{\max}-1}{\kappa_{\max}+1}.
\end{equation}
\end{theorem}
\begin{proof}
Only a small change in the proof of Theorem \ref{main_theorem} is needed. 
We first choose any vector $q'_0$, satisfying 
$\sin \angle_A(q_0,e_0)={(\kappa_{\max}-1)}/{(\kappa_{\max}+1)}.$
Then by Lemma \ref{neymeyr-like_lemma_complex} we obtain the complex Hermitian positive definite matrix $C_0$ and the complex scalar $\gamma$ such that $C_0e_0=\gamma q'_0$. Finally, we choose $q_0$ to be $\gamma q'_0$ and continue as in the proof of Theorem \ref{main_theorem}. 
The same modification is made in the choice of the vectors $q_k$ for $k\ge 1$ later in the proof.
\end{proof}

\section{Practical PCG algorithms}
\label{practical_pcg_sect}
In this section we briefly discuss two particular well-known PCG algorithms that are often used in practice. Our discussion here is motivated by and follows  
\citet[Remark 2.3]{golub_ye_MR1740397}. Suppose $A$, $b$, $x_0$, $r_0=b-Ax_0$, $\left\{B_k\right\}$ for $k=0,1,\ldots$ are given and consider Algorithm \ref{m2} where $\beta_k$ on line \ref{beta_line} is defined either by expression
\begin{equation}
\label{beta_standard}
\beta_k=\frac{\left(s_k,r_k\right)}{\left(s_{k-1},r_{k-1}\right)},
\end{equation}
or by expression
\begin{equation}
\label{beta_modif}
\beta_k=\frac{\left(s_k,r_k-r_{k-1}\right)}{\left(s_{k-1},r_{k-1}\right)}. 
\end{equation}
Formula (\ref{beta_standard}) is more often used in practice compared to 
(\ref{beta_modif}), since it can be implemented in such a way that 
does not require storing the extra vector $r_{k-1}$. 
\begin{algorithm}
\caption{}
\label{m2}
\begin{algorithmic}[1]
\FOR{$k=0,1,\ldots$}
\STATE $s_k=B_k^{-1} r_k$
\IF{$k=0$}
\STATE $p_0=s_0$
\ELSE
\STATE \label{beta_line} $p_k=s_k+\beta_k p_{k-1}$ \COMMENT{where $\beta_k$ is defined by either (\ref{beta_modif}) or (\ref{beta_standard}) for all iterations}
\ENDIF
\STATE $\D\alpha_k=\frac{\lb s_k,r_k\rb}{\lb p_k,Ap_k\rb}$
\STATE $\D x_{k+1}=x_k+\alpha_k p_k$
\STATE $\D r_{k+1}=r_k-\alpha_k Ap_k$
\ENDFOR
\end{algorithmic}
\end{algorithm}

If the preconditioner is SPD and fixed, it is well-known, e.g.,\ 
\citet[Remark 2.3]{golub_ye_MR1740397}, that $\D\lb s_k,r_{k-1}\rb=0$, so formula (\ref{beta_modif}) coincides with (\ref{beta_standard}) and  Algorithm \ref{m2} is described by (\ref{m}) with $m_k=\min\lb k,1\rb$. Of course, in this case the choice $m_k=\min\lb k,1\rb$ is enough to keep \emph{all} search directions $A$-orthogonal in exact arithmetic. 

Things become different when variable preconditioning is used.
It is well-known, e.g.,\ \citet[Remark 2.3]{golub_ye_MR1740397}
and \citet[Table 2]{notay_MR1797890}, 
that using formula (\ref{beta_standard}) for $\beta_k$ can significantly  slow down the convergence, 
and we provide our own numerical evidence of that in section \ref{num_exp_sect}. 
At the same time, comparing Lemma 3.2 with Lemma 2.1 from \citet{golub_ye_MR1740397},   	 
we can show, see \citet[v1]{kl06}, that 
Algorithm \ref{m2} with $\beta_k$ defined by (\ref{beta_modif}),
which is exactly \citet[IPCG Algorithm]{golub_ye_MR1740397},
is equivalent to the particular case of (\ref{m}), namely with $m_k=\min\lb k,1\rb$,  
and therefore is guaranteed by Theorem \ref{cr_th} to converge with at least the same speed as 
the PSD method. 

\section{Numerical experiments}
\label{num_exp_sect}
We first illustrate the main theoretical results of the paper
numerically for a model problem. 
We numerically investigate the influence of the 
choice for $\beta_k$ between
formulas (\ref{beta_standard}) and (\ref{beta_modif})
in Algorithm \ref{m2} and observe that (\ref{beta_modif}) leads to the theoretically 
predicted convergence rate, 
while (\ref{beta_standard}) may significantly  slow down the convergence.	 
Second, we test the convergence of inner-outer iteration schemes, where the 
inner iterations play the role of the variable preconditioning in the outer 
PCG iteration, and we illustrate our main conclusion that 
variable preconditioning may effectively 
reduce the convergence speed of the PCG method to the speed of the PSD method. 
Third, and last, we test the PSD and PCG methods with preconditioners 
of the same quality chosen randomly. We observe a surprising acceleration 
of the PCG method compared to the use of only one fixed preconditioner; 
at the same time, we show that the PSD method with random preconditioners 
works as well as the PCG method, which explains the PCG acceleration and 
again supports our main conclusion. 

\subsection{Numerical illustration of the main results} \label{num_exp_ssect1}
Here, we use the standard 3-point approximation of the 1-D Laplacian 
of the size 200 as the matrix $A$ of the system. To simulate the application of the variable preconditioner, we essentially repeat the steps described in the proof of Theorem \ref{main_theorem}, i.e.,\ we fix the condition number $\kappa\lb B_k^{-1}A\rb=2$ and on each iteration we generate a pseudo-random vector $s_k$, which is $A$-orthogonal to previous search directions 
and such that the $A$-angle between $s_k$ and $e_k$ satisfies 
$\sin\lb\angle_A\lb s_k,e_k\rb\rb=(\kappa-1)/(\kappa+1)$.

\begin{figure}
\centering
\includegraphics[scale=.4]{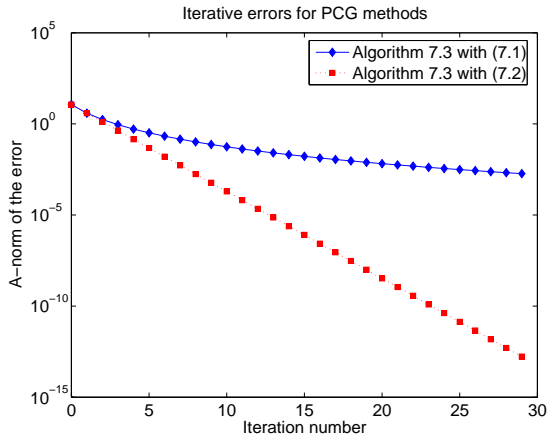}
\caption{\label{fig_alg_m2} Algorithm \ref{m2} with (\ref{beta_standard}) fails to provide the PSD convergence rate.}
\end{figure}
We summarize the numerical results of this subsection on Figure \ref{fig_alg_m2}, 
where the horizontal axis represents the number of iterations and the vertical axis represents the $A$-norm of the error. The iteration count actually starts from 1, so the $A$-norm of the error on the 0-th iteration $\vnorm{e_0}_A$ is just the $A$-norm of the initial error. 
The straight dotted (red in the colored print) line marked with squares on Figure \ref{fig_alg_m2} 
represents the PSD theoretical bound (\ref{classic_psd_bound}) 
and at the same time it perfectly coincides, which illustrates the statements of Theorem \ref{main_theorem}, with 
the change of the $A$-norm of the error in the case where the complete $A$-orthogonalization is performed, i.e.,\ $m_k=k$ in method (\ref{m}), 
as well as in the case where Algorithm \ref{m2} with $\beta_k$ defined by (\ref{beta_modif}) is used.  
The curved solid (blue) line marked with diamonds represents the convergence of Algorithm \ref{m2} with $\beta_k$ defined by (\ref{beta_standard}), which visibly performs much worse in this test compared to Algorithm \ref{m2} with (\ref{beta_modif}).  
The paper \citet[Sec. 5.2]{notay_MR1797890} contains analogous results 
comparing the change in the convergence rate using formulas 
(\ref{beta_standard}) and (\ref{beta_modif}), but 
it misses a comparison with the PSD method.  
To check our results of section \ref{complex_case_sect}, we repeat the tests
in the complex arithmetic. 
The figure generated is similar to Figure \ref{fig_alg_m2}, so we 
do not reproduce it here.  

\subsection{Inner-outer iterations as variable preconditioning} \label{num_exp_ssect2}
Inner-outer iterative schemes, where the 
inner iterations play the role of the variable preconditioner in the outer 
PCG iteration is a traditional example of variable preconditioning; see, e.g., 
\citet{golub_ye_MR1740397,notay_MR1797890}. 
Previously published tests analyze an approximation of some fixed preconditioner, $B_k \approx B$,  
different from $A$, by inner iterations, typically using the PCG method. The quality of the 
approximation is determined by the stopping criteria of the inner PCG method. 
A typical conclusion is that the performance of the outer PCG method
improves and it starts behaving like the PCG method with the fixed preconditioner $B$
when $B_k$ approximates $B$ more accurately by performing more inner iterations.  

\begin{figure}
\centering
\includegraphics[scale=.4]{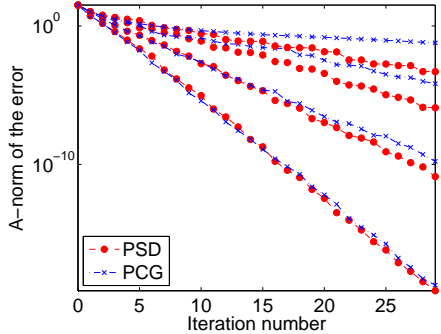}
\caption{\label{fig_inner_test} 
The PSD and PCG methods with preconditioning by inner CG 
with different stopping criteria $\eta= 0.2, 0.4, 0.6$ and $0.8$ (from the bottom to the top).}
\end{figure} 
The idea of our tests in this subsection is different: 
we approximate $B_k \approx B=A$. The specific setup is the following. 
We take a diagonal matrix $A$ with all integer entries from $1$ to $2000$,  
with the right-hand side zero and a random normally distributed zero mean
initial guess, and we do the same for the PSD and PCG methods.
For preconditioning on the $k$-th step, 
applied to the residual $r_k$, we run 
the standard CG method without preconditioning as inner iterations, 
using the zero initial approximation, and for the stopping criteria we 
compute the norm of the true residual at every inner iteration 
and iterate until it gets smaller than $\eta \| r_k\|$ for 
a given constant $\eta$. On Figure \ref{fig_inner_test}, 
we demonstrate the performance of the PSD and PCG methods for 
four values of $\eta= 0.2, 0.4, 0.6$ and $0.8$ (from the bottom to the top).
We observe that the PSD, displayed using dashed  (red in the colored print) lines marked with circles and 
PCG shown as dash-dot (blue) lines with x-marks  methods both converge with a similar rate, 
for each tested value of $\eta$.
We notice here that the PSD method is even a bit faster than the PCG method. 
This does not contradict our Corollary \ref{lo_csd}, since the preconditioners $B_k$ 
here are evidently different in the PSD and PCG methods even though they are 
constructed using the same principle. 
 
\subsection{Random vs. fixed preconditioning} \label{num_exp_ssect3}
In this subsection, we numerically investigate a situation where random preconditioners of a similar quality are used in the course of iterations. 
The system matrix is the standard 3-point finite-difference approximation of 
the one-dimensional Laplacian using $3000$ uniform mesh points 
and the Dirichlet boundary conditions.
We test the simplest multigrid preconditioning using two grids, 
where the number of coarse grid points is $600.$ 
The interpolation is linear, the restriction is the transpose of the interpolation, and the coarse-grid operator is defined by the Galerkin condition. The smoother is the Richardson iteration.
\begin{figure}
\centering
\includegraphics[scale=.35]{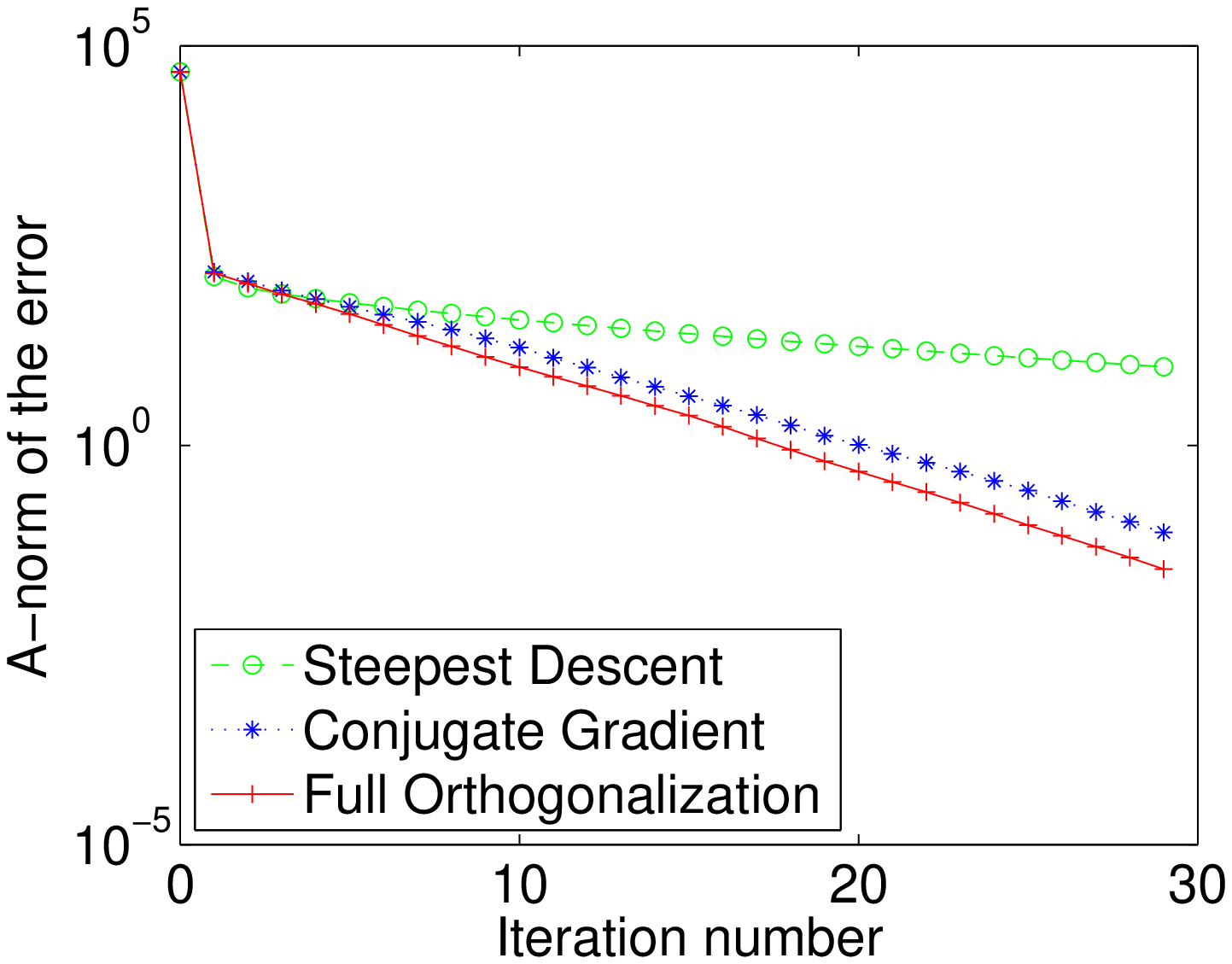}
\includegraphics[scale=.35]{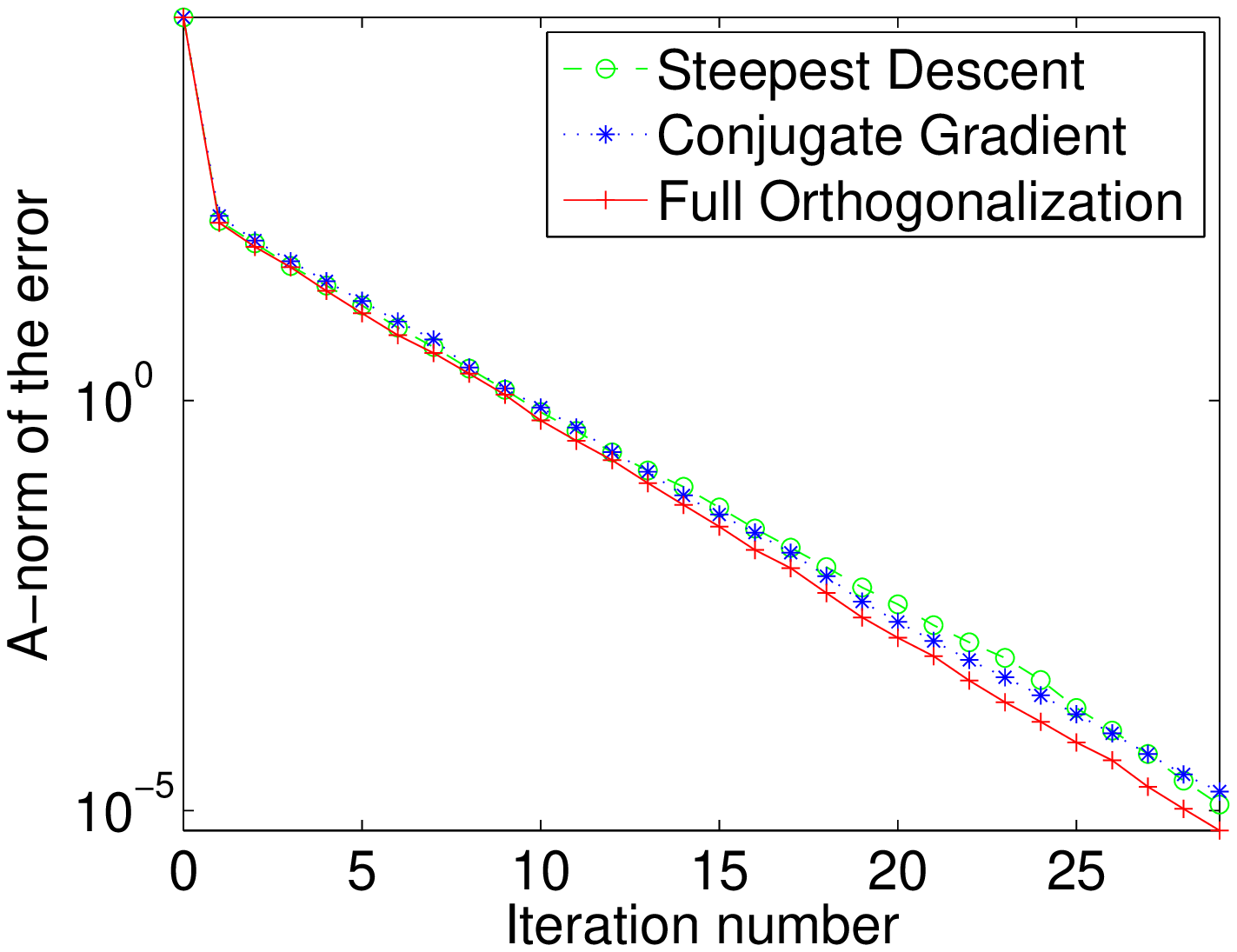}
\caption{\label{fig_random_multigrid} Two-grid preconditioning with fixed (left) and random (left) coarse grids.}
\end{figure} 
On Figure \ref{fig_random_multigrid} left, 
we once choose (pseudo-)randomly 600 coarse mesh points and   
build the fixed two-grid preconditioner, based on this choice. 
On Figure \ref{fig_random_multigrid} right,
we choose 600 new random coarse mesh points 
and rebuild the two-grid preconditioner on each iteration.
We note that in the algebraic multigrid the geometric information about the actual 
position of the coarse grid points is not available, so the random choice of the 
coarse grids may be an interesting alternative to traditional approaches. 

Figure \ref{fig_random_multigrid} displays the convergence history  
for the PSD (top), PCG (middle), and PCG with the full orthogonalization (bottom) 
with the same random initial guess using the fixed (left) and variable (right) 
two-grid preconditioners. 
On Figure \ref{fig_random_multigrid} left, for a fixed preconditioner, 
we observe the expected convergence behavior, with the PSD being noticeably the slowest, and the 
PCG with the full orthogonalization being slightly faster than the standard PCG. 
Figure \ref{fig_random_multigrid} right demonstrates that all three methods  
with the variable random preconditioner converge with essentially the same rate,  
which again illustrates the main result of the paper 
that the PCG method with variable preconditioning may just converge 
with the same speed as the PSD method. 

Figure \ref{fig_random_multigrid} reveals a surprising fact that 
the methods with random preconditioning converge twice as fast as 
the methods with fixed preconditioning! 
We highlight that Figure \ref{fig_random_multigrid} 
shows a typical case, not a random outlier, as we confirm 
by repeating the fixed preconditioner test in the left panel for \emph{every} 
random preconditioner used in the right panel of Figure \ref{fig_random_multigrid} and 
by running the tests multiple times with different seeds. 
Our informal explanation for the fast convergence of the PSD method with random preconditioning 
is based on Lemma \ref{sd_reduct_lemma} that provides 
the exact expression for the error reduction factor
as $\sin(\angle_A(e_k,B_k^{-1}Ae_k))$. It takes its largest value only 
if $e_k$ is one of specific linear combination of the eigenvectors of 
$B_k^{-1}Ae$ corresponding to the two extreme eigenvalues.
If $B_k$ is fixed, the error $e_k$ in the PSD method after several first iterations 
approaches these magic linear combinations, e.g.,\ \citet{MR0223071}, 
and the convergence rate reaches its upper bound. 
If $B_k$ changes randomly, as in our test, the average ``effective'' 
angle is smaller, i.e.,\ the convergence is faster. 
         
\section*{Conclusions}
We use geometric arguments to investigate the behavior 
of the PCG methods with variable preconditioning 
under a rather weak assumption that the quality of the preconditioner 
is fixed. Our main result is negative in its nature: we show that 
under this assumption the PCG method with variable preconditioning 
may converge as slow as the PSD method, moreover, as the PSD method 
with the slowest rate guaranteed by the classical 
convergence rate bound. In particular, 
that gives the negative answer, under our assumption,  
to the question asked in 
\citet[Sec. 6, Conclusion.]{golub_ye_MR1740397}
whether better bounds for the steepest descent reduction factor 
may exists for Algorithm \ref{m2} with (\ref{beta_modif}). 

Stronger assumptions on variable preconditioning, e.g.,\ such as 
made in \citet{golub_ye_MR1740397,notay_MR1797890}
that the variable preconditioners are all small 
perturbations of some fixed preconditioner,  
are necessary in order to hope to prove a convergence 
rate bound of the PCG method with variable preconditioning 
resembling the standard convergence 
rate bound of the PCG method with fixed preconditioning. 
Such stronger assumptions hold in many 
presently known real life applications 
of the PCG methods with variable preconditioning, but 
often require extra computational work, e.g.,\ 
more inner iterations in the inner-outer iterative methods.  



\bibliographystyle{plainnat}

\def\cprime{$'$} \def\cprime{$'$}
\begin{thebibliography}{15}
\providecommand{\natexlab}[1]{#1}
\providecommand{\url}[1]{\texttt{#1}}
\expandafter\ifx\csname urlstyle\endcsname\relax
  \providecommand{\doi}[1]{doi: #1}\else
  \providecommand{\doi}{doi: \begingroup \urlstyle{rm}\Url}\fi

\bibitem[Axelsson and Vassilevski(1994)]{MR1269945}
O.~Axelsson and P.~S. Vassilevski.
\newblock Variable-step multilevel preconditioning methods. {I}. {S}elfadjoint
  and positive definite elliptic problems.
\newblock \emph{Numer. Linear Algebra Appl.}, 1\penalty0 (1):\penalty0 75--101,
  1994.
\newblock ISSN 1070-5325.

\bibitem[Axelsson and Vassilevski(1991)]{axel_vass_MR1121697}
O.~Axelsson and P.~S. Vassilevski.
\newblock A black box generalized conjugate gradient solver with inner
  iterations and variable-step preconditioning.
\newblock \emph{SIAM J. Matrix Anal. Appl.}, 12\penalty0 (4):\penalty0
  625--644, 1991.
\newblock ISSN 0895-4798.

\bibitem[Axelsson(1994)]{axel_MR1276069}
Owe Axelsson.
\newblock \emph{Iterative solution methods}.
\newblock Cambridge University Press, Cambridge, 1994.
\newblock ISBN 0-521-44524-8.

\bibitem[D{\cprime}yakonov(1996)]{MR1396083}
Eugene~G. D{\cprime}yakonov.
\newblock \emph{Optimization in solving elliptic problems}.
\newblock CRC Press, Boca Raton, FL, 1996.
\newblock ISBN 0-8493-2872-1.
\newblock Translated from the 1989 Russian original, Translation edited and
  with a preface by Steve McCormick.

\bibitem[Forsythe(1968)]{MR0223071}
George~E. Forsythe.
\newblock On the asymptotic directions of the {$s$}-dimensional optimum
  gradient method.
\newblock \emph{Numer. Math.}, 11:\penalty0 57--76, 1968.
\newblock ISSN 0029-599X.

\bibitem[Golub and Ye(1999/00)]{golub_ye_MR1740397}
Gene~H. Golub and Qiang Ye.
\newblock Inexact preconditioned conjugate gradient method with inner-outer
  iteration.
\newblock \emph{SIAM J. Sci. Comput.}, 21\penalty0 (4):\penalty0 1305--1320
  (electronic), 1999/00.
\newblock ISSN 1064-8275.

\bibitem[Kantorovi{\v{c}}(1947)]{MR0023126}
L.~V. Kantorovi{\v{c}}.
\newblock On the method of steepest descent.
\newblock \emph{Doklady Akad. Nauk SSSR (N. S.)}, 56:\penalty0 233--236, 1947.

\bibitem[Kantorovich and Akilov(1964)]{Kantorovich-1964}
L.~Kantorovich and G.~P. Akilov.
\newblock \emph{Functional Analysis in Normed Spaces}.
\newblock Pergamon, NY, 1964.

\bibitem[Knyazev and Lashuk(2006-2007)]{kl06}
Andrew~V. Knyazev and Ilya Lashuk.
\newblock Steepest descent and conjugate gradient methods with variable
  preconditioning.
\newblock Electronic. math.NA/0605767, arXiv.org,
  http://arxiv.org/abs/math/0605767, 2006-2007.
\newblock Revised v3.

\bibitem[Neymeyr(2001)]{neymeyr_MR1804114}
Klaus Neymeyr.
\newblock A geometric theory for preconditioned inverse iteration. {I}.
  {E}xtrema of the {R}ayleigh quotient.
\newblock \emph{Linear Algebra Appl.}, 322\penalty0 (1-3):\penalty0 61--85,
  2001.
\newblock ISSN 0024-3795.

\bibitem[Notay(2000)]{notay_MR1797890}
Yvan Notay.
\newblock Flexible conjugate gradients.
\newblock \emph{SIAM J. Sci. Comput.}, 22\penalty0 (4):\penalty0 1444--1460
  (electronic), 2000.
\newblock ISSN 1064-8275.

\bibitem[Simoncini and Szyld(2002)]{MR1974182}
Valeria Simoncini and Daniel~B. Szyld.
\newblock Flexible inner-outer {K}rylov subspace methods.
\newblock \emph{SIAM J. Numer. Anal.}, 40\penalty0 (6):\penalty0 2219--2239
  (electronic) (2003), 2002.
\newblock ISSN 0036-1429.

\bibitem[Simoncini and Szyld(2003)]{MR2058070}
Valeria Simoncini and Daniel~B. Szyld.
\newblock Theory of inexact {K}rylov subspace methods and applications to
  scientific computing.
\newblock \emph{SIAM J. Sci. Comput.}, 25\penalty0 (2):\penalty0 454--477
  (electronic), 2003.
\newblock ISSN 1064-8275.

\bibitem[Simoncini and Szyld(2005)]{MR2179897}
Valeria Simoncini and Daniel~B. Szyld.
\newblock On the occurrence of superlinear convergence of exact and inexact
  {K}rylov subspace methods.
\newblock \emph{SIAM Rev.}, 47\penalty0 (2):\penalty0 247--272 (electronic),
  2005.
\newblock ISSN 0036-1445.

\bibitem[Szyld and Vogel(2001)]{MR1861254}
Daniel~B. Szyld and Judith~A. Vogel.
\newblock F{QMR}: a flexible quasi-minimal residual method with inexact
  preconditioning.
\newblock \emph{SIAM J. Sci. Comput.}, 23\penalty0 (2):\penalty0 363--380
  (electronic), 2001.
\newblock ISSN 1064-8275.
\newblock Copper Mountain Conference (2000).

\end{thebibliography}

\def\cprime{$'$} \def\cprime{$'$}

\end{document}